\documentclass[fontsize=12pt,parskip=full]{scrartcl}
\usepackage{amsmath,amssymb,amsfonts,amsthm,mathabx}
\usepackage[dvips]{graphics}
\usepackage{enumerate}
\usepackage{mathrsfs}
\usepackage{subfigure}
\usepackage{color}

\usepackage[T1]{fontenc}

\DeclareSymbolFont{AMSb}{U}{msb}{m}{n}

\def\P{\mathbb P}

\def\R{\mathbb {R}}

\def\F{\mathcal {F}}

\def\Ecal{\mathcal {E}}
\def\E{\mathbb {E}}

\def\1{\,{\makebox[0pt][c]{\normalfont    1}
\makebox[2.5pt][c]{\raisebox{3.5pt}{\tiny {$\|$}}}
\makebox[-2.5pt][c]{\raisebox{1.7pt}{\tiny {$\|$}}}
\makebox[2.5pt][c]{} }}

\def\scirc{\mathbin{\raise.15ex\hbox{\scriptsize$\circ$}}}

\newtheorem{thm}{Theorem}

\theoremstyle{definition}
\newtheorem{defn}[thm]{Definition}

\title{Minimal forward random point attractors need not exist}

\author{Michael Scheutzow%
  \thanks{Institut f\"ur Mathematik, MA 7-5, Fakult\"at II, 
    Technische Universit\"at Berlin, 
    Stra\ss e des 17.~Juni 136, 10623 Berlin, FRG;  \ 
    \small\texttt{ms{\scriptsize @}math.tu-berlin.de}}}

\date{\today}

\begin{document}\maketitle

\begin{abstract}\noindent
  It is well-known that random attractors of a random dynamical system are generally not unique.
  It was shown in \cite{CS18}  that if there exist more than one pullback or weak random attractor which attracts a given family of
  (possibly random) sets, then there exists a minimal (in the sense of smallest) one. This statement does not hold for forward
  random attractors. The same paper contains an example of a random dynamical system and a deterministic family of sets
  which has more than one forward attractor which attracts the given family but no minimal one. The question whether one can find an
  example which has multiple forward {\em point} attractors but no minimal one remained open. Here we provide such an example.

\par\medskip

  \ Secondary\, 
37B25 \ 
37H99 \ 
37L55 \ 
\end{abstract}


\section{Introduction}
For random dynamical systems on a Polish (i.e.~complete separable metric) space $E$ random attractors have been studied extensively during the past 25 years.
There are several different concepts of a random attractor depending on the family of sets which are attracted as well as the mode of attraction: weak (i.e.~attraction
in probability), pullback and forward. The most common considered families of attracted sets are the family of all deterministic compact (or bounded) sets
as well as the family of all deterministic singletons. The former are called {\em set} attractors and the latter {\em point} attractors. While
(weak, pullback and forward) set attractors are known to be unique, this is not true for point attractors and the question arises if there always exists a smallest point
attractor (in case we know that at least one such attractor exists). For pullback and weak point attractors a positive answer was provided in \cite{CS18}.
The same paper also contains an example of a random dynamical system and a family of deterministic sets for which there exist many forward attractors but not a smallest one.
The question whether a minimal forward  {\em point} attractor always exist (in case there exists at least one such attractor) remained open. Here we provide an example which
shows that the answer is negative. The example is an appropriate modification of the example given in \cite{CS18}.

Before we present the example, we recall the definition of a forward point attractor. For the definition of a random dynamical system and other kind of attractors, we
refer the reader to \cite{CS18}. For simplicity we only consider the continuous time case.

\begin{defn} Let $(\Omega,\F,\P,\vartheta,\varphi)$ be a continuous time random dynamical system taking values in the Polish space $E$ with complete metric $d$.
  Denote the Borel $\sigma$-algebra on $E$ by $\Ecal$.

Then a set $A\subset E\times\Omega$ is a \emph{(random) forward point attractor} if
\begin{enumerate}[(i)]
\item $A$ is a compact random set (i.e.~$A \in \Ecal \otimes  \F$ and all sections $A(\omega)$ are compact)
\item $A$ is strictly $\varphi$-invariant, i.e.\ 
  \begin{equation*}
    \varphi_t(\omega)A(\omega)=A(\vartheta_t\omega)
  \end{equation*}
   $\P$-almost surely for every $t\geq0$ 
\item $A$ attracts points, i.e. 
  \begin{equation*} 
    \lim_{t\to\infty}
    d\bigl(\varphi_t(\omega)x,A(\vartheta_t\omega)\bigr)=0\qquad\P\mbox{-a.s.},
  \end{equation*}
  for every $x \in E$. 
\end{enumerate}
\end{defn}

\section{The example}

In this section, we provide an example of an RDS which has more than one  forward point attractor but not a smallest one.   

Consider a stationary Ornstein-Uhlenbeck process $Z$, i.e.~a real-valued centered Gaussian process defined on $\R$ with covariance $\E (Z_tZ_s)=\frac 12\exp\{-|t-s|\}$. 
We define $Z$ on the 
canonical space $\Omega=C(\R,\R)$ of continuous functions from $\R$ to $\R$ together with the usual shift $\vartheta$ and equipped with the law of $Z$.
Then $Z_t(\omega)=Z_0(\vartheta_t\omega)$. Let $\tilde d$ be the Euclidean metric on $\R$, let $E:=\R \times [0,1]$ be equipped with the Euclidean metric $d$, and define
$\Gamma:=\{-3,-1,1,3\}$.

For $(x,y) \in E$ let $h(t,x,y)$, $t \ge 0$ be the unique solution of the ode 
$$
\dot u=u-1-\tilde d(x,\Gamma)\wedge 1
$$
with initial condition $h(0,x,y)=y \in [0,1]$.


Next, we define, for $x \in \R$, $y \in [0,1]$, and $t \ge 0$
\begin{align*}
\varphi_t(\omega)(x,y):=
\begin{cases}
(x+Z_t(\omega)-Z_0(\omega),\,h(t,x-Z_0(\omega),y)), &\mbox{if } t \le \tau(x-Z_0(\omega),y)\\   
\big(\exp\{ \tau(x-Z_0(\omega),y)         -t  \}\big(  x-Z_0(\omega)\big)+Z_t(\omega),0\big), &\mbox{if } t \ge  \tau(x-Z_0(\omega),y),
\quad 
\end{cases}
\end{align*}
where
$$
\tau(x,y):=\inf\{s \ge 0: h(s,x,y)=0\}.
$$
It is straightforward to check that $\varphi$ defines a continuous random dynamical system on $\Omega$. Note that at time $\tau(x-Z_0(\omega),y)$, the process $\varphi_t(x,y)$ starting in $(x,y)$ hits the $x$-axis and then moves on the $x$-axis approaching the process
$Z_t(\omega)$ with deterministic exponential speed. The only initial points for which $\tau(x-Z_0(\omega),y)=\infty$ i.e.~for which the trajectories will never hit the $x$-axis are
the four points in $G(\omega):=\{(Z_0(\omega)+\gamma,1);\,\gamma \in \Gamma\}$.
All trajectories starting outside this set will converge to $(Z_t(\omega),0)$ with respect to  the Euclidean metric $d$ on $E$. Since for any deterministic $(x,y)\in E$ we have
$\P((x,y)\in G(\omega))=0$ we see that $A(\omega):=\{(Z_0(\omega),0)\}$ is a random forward point attractor with respect to $d$. 
This is no longer true if we change the metric on $E$ in the following way (without changing the topology of $E$):
$$
\rho\big((x,y),(\tilde x,\tilde y)\big):= |\tilde y -y|+|H (\tilde x)-H(x)|,
$$   
where $H$ is strictly increasing, odd, continuous such that $H(x)=\exp\{\exp\{\exp (x)\}\}$ for large $x$ (the fact that this metric  works can be checked 
by using the fact that the running maximum of a stationary Ornstein Uhlenbeck process up to time $t$ is of the order $\sqrt{\log t}$). Then
$\limsup_{t \to \infty} \rho(\varphi_t(x,y),(Z_t(\omega),0))=\infty$ almost surely for every fixed $(x,y) \in E$ and therefore, in particular, $A(\omega):=\{(Z(\omega),0)\}$ is not
a forward point attractor with respect to $\rho$ (it does not even attract any point). There are however several forward 
point attractors with respect to the metric $\rho$, for example 
$$
A_\gamma(\omega):= \big([Z_0(\omega)-\gamma,Z_0(\omega)+\gamma]\times \{0\}\big)\bigcup  \big(\{Z_0(\omega)-\gamma\}\times [0,1]\big) \bigcup  \big( \{Z_0(\omega)+\gamma\}\times [0,1]\big)
$$ 
for any $\gamma \in \{1,3\}$ (note that these two sets are strictly invariant!). These two random sets do not only attract but even ``swallow'' (or absorb) the trajectory starting at
$(x,y)$ for every $(x,y)\in E$ almost surely. Assume that there is a smallest forward point attractor $\hat A(\omega)$ with respect to $\rho$.
Then $\hat A(\omega)$ has to be contained in the intersection 
$A_1(\omega)\bigcap A_3(\omega)$ which is a subset of $\R \times \{0\}$. It is clear however that the set $\{(Z_0(\omega),0)\}$ is the only strictly invariant compact subset of $\R \times \{0\}$ 
and we already saw that this is not a forward attractor, contradicting our assumption. Hence, this RDS does not have a smallest forward point attractor.


\begin{thebibliography}{10}%
  \bibliographystyle{abbrv}
 \bibitem{CS18}
  H.~Crauel and M.~Scheutzow, 
  Minimal random attractors, 
  \emph{J. Differential Equations} 265 (2018) 702-718   
\end{thebibliography}
\end{document}